# Complexity and asymptotic stability in the process of biochemical substance exchange in a coupled ring of cells

D.T. Mihailović[1], V. Kostić[2], I. Balaž[3] and Lj. Cvetković[2]

[1] Faculty of Agriculture, Division of Meteorology and Biophysics, University of Novi Sad, Dositeja Obradovica Sq. 8, 21000 Novi Sad, Serbia
[2] Faculty of Sciences, Department of Mathematics and Informatics, University of Novi Sad, Dositeja Obradovica Sq. 3, 21000 Novi Sad, Serbia
[3] Faculty of Sciences, Department of Physics, University of Novi Sad, Dositeja Obradovica Sq. 3, 21000 Novi Sad, Serbia

**Abstract**

We have considered the complexity and asymptotic stability in the process of biochemical substance exchange in a coupled ring of cells. We have used coupled maps to model this process. It includes the coupling parameter, cell affinity and environmental factor as master parameters of the model. We have introduced: (i) the Kolmogorov complexity spectrum and (ii) the Kolmogorov complexity spectrum highest value to analyze the dynamics of two cell model. The asymptotic stability of this dynamical system using an eigenvalue-based method has been considered. Using these complexity measures we have noticed an "island" of low complexity in the space of the master parameters for the weak coupling. We have explored how stability of the equilibrium of the biochemical substance exchange in a multi-cell system ($N=100$) is influenced by the changes in the master parameters of the model for the weak and strong coupling. We have found that in highly chaotic conditions there exists space of master parameters for which the process of biochemical substance exchange in a coupled ring of cells is stable.

**Keywords:** complexity, stability, coupling, substance exchange

## 1.Introduction

Synchronization between different processes is one of the prerequisites for achieving functionality in the cell. Moreover, question how synchronization influence robustness, adaptability and evolvability of multi-cell systems is one of the most important topics of contemporary biophysics of complex systems [1-6]. The exchange of biophysical substances among the components of multi-cell systems is driven by a range of both intra- and inter-cellular factors. Recently, significant efforts have been made toward better understanding of the stability of synchronization in biological systems and its influence on the functional robustness of multi-cell systems [2-5, 7-9]. In those papers the authors considered cells as uniform particles, without internal structure and without the ability to change their behavior. However, in natural conditions, cells usually shift between several growth and development

□ Corresponding author (D.T. Mihailović)

*Telephone*: +381216350552

*E-mail address*:guto@polj.uns.ac.rs

phases. For example, bacterial cells spend most of the time in the stationary phase which is characterized by a decreased growth rate, slowdown of all metabolic processes and increase in resistance to several stress conditions [10-13]. In our previous work we investigated how change of affinity in the coupling affects synchronization of substance exchange between cells [6]. However, it still remains open question how the complexity and stability of the substance exchange process are affected by the changes in parameters that represent the influence of the environment, cell coupling and cell affinity.

The issue of system complexity in many disciplines, ranging from cognitive science to evolutionary and developmental biology and particle astrophysics [14-18 and references herein], has been touching the scientific community intensively during the last three decades. An important part in studying complex systems is analysis of symbolic sequences. Namely, it is believed that most systems whose complexity we would like to estimate can be reduced to them. However, according to Adami and Cerf [19], the idea that the regularity of the string alone is connected to its complexity, is meaningless if this analysis is done in the absence of an environment within which the string is to be interpreted. Thus, the complexity of a string representing the complexity of the mentioned systems can be determined only by analyzing its correlation with a system environment [19]. In this paper it will be a sequence of a system's process, as an indication of its complexity. The stability of the complex systems is also an important issue that drawn considerable attention in analysis of biological and infection-related systems [20-22, and herein references]. Here, we will analyze asymptotic stability of the dynamical system. Namely, if solutions that are sufficiently close to the equilibrium point remains sufficiently close to that point forever, then the system is said to be Lyapunov stable. If the system is Lyapunov stable and solutions that start sufficiently close to the equilibrium point eventually converge to the equilibrium, then the system is said to be asymptotically stable. The stability analysis for dynamics of different complex systems [23-26, among others] can be carried out by using different methods, but the eigenvalue-based methods are commonly used for the asymptotic stability.

In this paper, in order to provide a better insight in these problems, we will model the process of the biochemical substance exchange using a map, which includes a formalization of the cell affinity in the interaction with other cells. First, for two-cell system, we consider the complexity of this process by introducing the Kolmogorov complexity that has been applied on the modeled process, and then we discuss its asymptotic stability, through the theory of discrete time-invariant dynamical systems. And second, for the multi-cell system presented by a ring of coupled cells, we will provide set of conditions that describe how stability of the equilibrium of the biochemical substance exchange in time is influenced by the parameters of the model.

**2. Model background**

*2.1. Short description of the inter-cellular biochemical substance exchange model*

In order to explore asymptotic stability in the model of dynamics of biochemical substance exchange in a coupled ring of cells, first we shortly describe an inter-cellular exchange model (Figure 1). Signaling molecules are ones which are deliberately extracted by the cell, and can affect behavior of other cells of the same or different type by means of active uptake and subsequent physiological response. In this context, we can define two environments: (i) intra-cellular environment (inside the cell) and (ii) inter-cellular environment, which refers to the environment between cells. Since active uptake is one of the milestones of the process, a very important factor in establishing communication is set of active receptors and transporters in cellular membrane. Also, receptors are very important source of perturbations of the communication process due to protein disorder and intrinsic noise. Another important factor is inter-cellular environment which could interfere with the process of exchange. It includes: distance between cells, mechanical and dynamical properties of the fluid which serves as a channel for exchange and various abiotic and biotic factors influencing physiology of the involved cells. Finally, in order to define exchange process as communication, received molecules should induce physiological response. Therefore, concentration of signaling molecules inside of the cell can serve as an indicator of dynamics of the whole process of communication. Additionally, the influence of affinity in functioning of living systems is also an important issue. In this case, it can be divided into following aspects: (A1) affinity of genetic regulators towards arriving signals which determine intensity of cellular response and (A2) affinity for uptake of signaling molecules [27]. First aspect is genetically determined and therefore species specific. Second aspect is more complex and is influenced by: affinity of receptors to binding specific signaling molecules, number of active receptors and their conformational fluctuations (protein disorder).

As it is obvious from the empirical description, we can infer successfulness of the communication process by monitoring: (i) number of signaling molecules, both inside and outside of the cell and (ii) their mutual influence. Concentration of those molecules in inter-cellular environment depends on various environmental factors and taken alone often can indicate more about state of the environment then about the communication itself. Therefore, we choose the concentration of signaling molecules inside of the cell to be the main indicator of the communication process. With this choice, parameters of the system are: (i) cell affinity $p$ (with this parameter we formalize an intrinsic property of the cell [27,28]) by which cells perform uptake of signaling molecules (A2), that depends on number and state of appropriate receptors; (ii) concentration $c$ of signaling molecules in inter-cellular environment within the radius of interaction; (iii) intensity of cellular response (A1) for two cells $x_n$ and $y_n$ in discrete time step $n$ and (iv) influence of other environment factors $r$ which can interfere with the process of communication. Here, we postulate that parameter $r$ can be taken collectively for intra- and inter-cellular environment, inside of the one variable, indicating overall disposition of the environment to the communication process.

Then, the evolution of the concentration in two-cell system through time ($x_n$, $y_n$), $n = 0,1,2,...$, can be expressed as

$$x_{n+1} = (1-c)\varphi(x_n) + c\psi(y_n), \tag{1a}$$

$$y_{n+1} = (1-c)\varphi(y_n) + c\psi(x_n). \tag{1b}$$

Here, in a time step $n$, the map $\psi:(0,1) \to (0,1)$ represents the flow of materials from cell to cell. Maps $\psi(x_n)$ and $\psi(y_n)$ can be approximated by a power map, such that $\psi(x_n) \sim (x_n)^p$ and $\psi(y_n) \sim (y_n)^q$, where the affinity to uptake signaling molecules is indicated by parameters $p$ and $q$, respectively. Let us note that we require that sum of all the affinities should equal 1, or in the case of two cells $p+q=1$. In that way the interaction is expressed as a nonlinear coupling between two cells. On the other hand, the map $\varphi:(0,1) \to (0,1)$ models the dynamics of intra-cellular behavior, and it can be expressed as a logistic map $\varphi(x) = r\,x\,(1-x)$, where $r$ is a logistic parameter, $0 < r \leq 4$, [6, 29, 30], that indicates overall disposition of the environment to the communication process. Finally, $c$ represents coupling of two factors:

concentration of signaling molecules in intra-cellular environment and intensity of response they can provoke. This form is taken because the effect of the same intra-cellular concentration of signaling molecules can vary greatly with variation of affinity of genetic regulators for that signal, which is further reflected on the ability to synchronize with other cells.

Therefore, we consider system of difference equations in the vector form

$$\mathbf{v}_{n+1} = \mathbf{F}(\mathbf{v}_n) \equiv \mathbf{L}(\mathbf{v}_n) + \mathbf{P}(\mathbf{v}_n), \qquad (2)$$

with notation

$$\mathbf{L}(\mathbf{v}_n) = \begin{bmatrix} (1-c)rx_n(1-x_n) & (1-c)ry_n(1-y_n)) \end{bmatrix}^T, \qquad \mathbf{P}(\mathbf{v}_n) = \begin{bmatrix} c\, y_n^p & c\, x_n^{1-p} \end{bmatrix}^T, \quad (3)$$

where $\mathbf{v}_n = \begin{bmatrix} x_n & y_n \end{bmatrix}^T$ is a vector representing concentration of signaling molecules inside of the cells. It is important to notice that $\mathbf{P}(\mathbf{v}_n)$ denotes stimulating coupling influence of members of the system which is here restricted only to positive numbers in the interval (0,1). The starting point $\mathbf{v}_n$ is always chosen so that $\mathbf{v}_0 \in (0,1) \times (0,1)$. And, since we have a trivial fixed point $F(0) = 0$, in order to ensure that zero is not at the same time the point of attraction, we defined $p \in (0,1)$ as an exponent. Since $c$ influences both, rate of intra-cellular synthesis of signaling molecules and the synchronization of signaling processes between two cells, parameter $c$ is taken to be a part of both $\mathbf{L}(\mathbf{v}_n)$ and $\mathbf{P}(\mathbf{v}_n)$. However, relative ratio of these two influences depends on current model setting. Therefore, in this paper, we will first consider two-cell model and the case where both cell are strongly influenced by intra-cellular concentration of signals, while they can provoke relatively smaller response

$$x_{n+1} = (1-c)\, r\, x_n(1-x_n) + c\, y_n^p, \qquad (4a)$$

$$y_{n+1} = (1-c)\, r\, y_n(1-y_n) + c\, x_n^{1-p}. \qquad (4b)$$

*2.2. Complexity of the inter-cellular biochemical substance exchange model*

The notion of complexity is used in different ways, but overall consensus of its definition still does not exist. According to Grassberger [31] it can be intuitively accepted as something that is placed between uniformity and total randomness. According to Adami [32] the biophysical complexity of a sequence "refers to the amount of information that is stored in that sequence about a particular environment". This complexity should not be confused with mathematical (Kolmogorov) complexity; it is a distinct mathematical complexity, which only deals either with the intrinsic regularity or irregularity of a sequence in this case. Namely, for any two strings $x, y \in P^*$ ($P^*$ is the set of all finite binary strings), the Kolmogorov complexity of given $x$ is $K(x|y) = \min_p \{|p| : U(p, y) = x\}$, where $U(p, y)$ is the output of the program $p$ with auxiliary input $y$ when it is run in the machine $U$. A good introduction to the Kolmogorov complexity $K^C$ can be found in Cover and Thomas [33] and with a comprehensive description in Li and Vitanyi [34]. On the basis of Kolmogorov's idea, Lempel and Ziv [35] developed an algorithm, which is used in assessing the randomness of finite sequences as a measure of its disorder. A short description of this algorithm is given in the Appendix A. Let us note that the Kolmogorov complexity is not able to distinguish between time series which have different amplitude oscillations but very similar random components.

We analyze the complexity of the map (Eq. (4)), using the Kolmogorov complexity spectrum highest value, $K_m^C$ of the concentration time series $\{x_i\}$, $i = 1, 2, 3, 4, ..., N$, which we introduce in the following way. First, we transform this time series into a finite symbol string by comparison with series of thresholds $\{x_{t,i}\}$, $i = 1, 2, 3, 4, ..., N$, where each element is equal to the corresponding element in the considered time series $\{x_i\}$. The original time series samples are converted into a 0-1 sequences $\{S_i^{(k)}\}$, $k = 1, 2, 3, 4, ..., N$ defined by comparison with a threshold $x_{t,k}$,

$$S_i^{(k)} = \begin{cases} 0 & x_i < x_{t,k} \\ 1 & x_i \geq x_{t,k} \end{cases}. \qquad (5)$$

After we apply the Lempel-Ziv algorithm on each element of the time series $\{S_i^{(k)}\}$ we get the Kolmogorov complexity spectrum $\{K_i^C\}$, $i = 1, 2, 3, 4, ..., N$. We introduce this spectrum in order to explore the range of amplitudes in the time series representing a process of inter-cellular biochemical substance exchange, for which that process has the highest complexity. The highest value $K_m^C$ in this series, i.e. $K_m^C = \max\{K_i^C\}$, we call the Kolmogorov complexity spectrum highest value. Let us note that $K^C$ means a single value, which is computed using the threshold as explained in Appendix A.

Figure 2 shows the change of $K_m^C$ against cell affinity $p$ and logistic parameter $r$, with: (i) $r$ and $p$ taking values in the interval (3.6 - 4.0) and (0 - 1.0), respectively, with an increment of 0.005 and (ii) initial condition $(x_0, y_0) = (0.3, 0.5)$. For each $x$, 5000 iterations of the map (Eq. 4) are applied, and the first 1000 steps are abandoned. Figure 3a show that for $c = 0.02$ (weak coupling), $r > 3.7$ and whole range of the cell affinity $p$, $K_m^C$ takes high values corresponding to the high complexity in the simulation of process of the biochemical substance exchange between two cells. Only exception is the "island" around the point (3.95, 0.5) in ($r, p$) phase space, which can be attributed to the nonlinear features of the map (Eqs. 4a-4b). In this case, the high complexity indicates the presence of randomness in time series of concentration $x$ in the cell. In contrast to this, for higher values of coupling (starting from $c = 0.2$) the process of the biochemical substance exchange (Fig. 2a) is much less complex indicating absence of stochastic processes in comparison with the weak coupling (Fig. 2b).

Figure 3 depicts $K^C$ of the process of the biochemical substance exchange between two cells (Eq. (4a)-(4b)) as a function of normalized concentration in cell, $X_i$ for two values of the coupling parameter $c$ (0.02 and 0.2). We have obtained the normalized values of this series by the transformation $X_i = (x_i - x_{min}) / (x_{max} - x_{min})$, where $\{x_i\}$ is a time series of concentration in the cell. From this figure it is seen that the process of the biochemical substance exchanges, for weak coupling $c = 0.02$, happens for any concentration in the cells. However, for $c = 0.2$ substance exchange is located in the regions of $X_i$ - (0.4, 0.6) and (0.9, 1.0). Let us note that Hu et al. [36] derived analytic expression for $C_k$ (notation in Appendix A) in the Kolmogorov complexity, for regular and random sequences. In addition they showed that for the shorter length of the time series, the $C_k$ value is larger and correspondingly the complexity for a random sequence can become considerably larger than 1.

*2.3. Stability of the inter-cellular biochemical substance exchange model*

In this subsection we address the behavior of the system (2), and estimate whether this dynamical system can achieve the asymptotic stability at the provided equilibrium point $\tilde{\mathbf{v}} = \mathbf{F}(\tilde{\mathbf{v}})$ for two cells in dependence of coupling parameter $c$ and different values of system parameters $r$ and $p$. Since, in general, our map $\mathbf{F}:(0,1)\times(0,1) \to (0,1)\times(0,1)$ is a nonlinear, we will discuss the asymptotic stability behavior of the system (2) at the given equilibrium by the linearization technique. Namely, it is well known that the asymptotic stability of a dynamical system at the equilibrium point is described by the spectral behavior of its linearization around that equilibrium, [37]. Therefore, we will compute the Jacobian of the map $\mathbf{F}$ at the point $\tilde{\mathbf{v}} = [\tilde{x}\ \ \tilde{y}]^T$ as

$$J_F(\tilde{\mathbf{v}}) := D\mathbf{F}|_{\mathbf{v}=\tilde{\mathbf{v}}} = \begin{bmatrix} (1-c)r(1-2\tilde{x}) & cp(\tilde{y})^{p-1} \\ c(1-p)(\tilde{x})^{-p} & (1-c)r(1-2\tilde{y}) \end{bmatrix}, \quad (6)$$

with eigenvalues

$$\lambda_1 = (1-c)r(1-\alpha(\tilde{x},\tilde{y})) + \sqrt{(1-c)^2 r^2 (1-\alpha(\tilde{x},\tilde{y}))^2 + c^2(1-c)^2 \frac{\gamma(\tilde{x},\tilde{y})}{\beta(\tilde{x},\tilde{y})}}, \quad (7a)$$

$$\lambda_2 = (1-c)r(1-\alpha(\tilde{x},\tilde{y})) - \sqrt{(1-c)^2 r^2 (1-\alpha(\tilde{x},\tilde{y}))^2 + c^2(1-c)^2 \frac{\gamma(\tilde{x},\tilde{y})}{\beta(\tilde{x},\tilde{y})}}, \quad (7b)$$

where: $\alpha(\tilde{x},\tilde{y}) := \tilde{x} + \tilde{y}$ is the total concentration in both cells ($\tilde{x}$ and $\tilde{y}$) at equilibrium state, $\beta(\tilde{x},\tilde{y}) := (\tilde{x}^p \tilde{y}^{1-p})/(p(1-p))$ is the weighted generalized geometric mean of concentrations of substance between two cells and $\gamma(\tilde{x},\tilde{y}) := (1-2\tilde{x})(1-2\tilde{y})$ whose sign indicates in which two quadrants, of the domain corresponding to the axis through the point (0.5, 0.5), the provided equilibrium state belongs to.

Since the matrix in (6) has only two eigenvalues, using Theorem from Appendix B, then the discrete nonlinear dynamical system (2) at the equilibrium point $\tilde{\mathbf{v}} = [\tilde{x}\ \ \tilde{y}]^T$ is

- asymptotically stable for

$$\max\{|\lambda_1|,|\lambda_2|\} < 1, \quad (8a)$$

and

- not asymptotically stable if:

$$\min\{|\lambda_1|, |\lambda_2|\} > 1, \tag{8b}$$

while otherwise we cannot get the answer through analysis of eigenvalues of the Jacobian. Therefore, (8a) and (8b) completely express the stability and instability that is obtained through the linearization. If these conditions are not satisfied, stability purely depends on nonlinearity of the map $\mathbf{F}$ given by (2), around the point of equilibrium $\tilde{\mathbf{v}}$. In this case the regions of stability and instability can be further explored by eigenvalues of Hessian, i.e. using second order of derivatives.

When for a given set of parameters expressions (7a) and (7b) are computed, then regions of the domain, for which equilibrium state of the system (2) is asymptotically either stable or unstable, can be plotted. In order to better grasp the meaning of stability/instability conditions given by (8a) and (8b) in Fig. 4 (week coupling $c = 0.02$) and Fig. 5 (strong coupling $c = 0.6$), we have plotted regions of the asymptotic stability and instability in the domain of the map $\mathbf{F}$ for the following values of parameters: $p = 0.25, 0.5, 0.75$ and $r = 1, 2, 3, 4$. In each plot of these figures, the values of the equilibrium concentration of substance in the first cell ($\tilde{x}$) are given on the horizontal axis, while the values of the second one ($\tilde{y}$) are given on the vertical axis. The dark gray areas of the domain $(0,1) \times (0,1)$, represent asymptotic stability, while light gray areas indicate regions of the instability. Note that white area corresponds to the region where stability/instability is due to purely nonlinear effects of the map $\mathbf{F}$ given by Eqs. (2).

In Fig. 4 for logistic parameter $r$ we use the broad range of its values from 1 to 4, although chaotic fluctuations in this equation occur after $r = 3.57$. From this figure it is seen that the asymptotic instability occurs over the whole domain (i.e., dark gray color completely covers the square section) for $r = 1$ and all values of $p$ (0.25, 0.50, 0.75). With increase of $r$ there exist three sections corresponding to asymptotic stability (dark gray), instability (light gray) and stability purely depending on nonlinearity of the map $\mathbf{F}$ (white). Note that the stability square section becomes smaller with increase of $r$. This is expected for the weak coupling since we approach to region with the chaotic fluctuations. Now we can explain the occurrence of a "island" around the point (3.95,0.5) in ($r, p$) phase space in Fig. 2, which depicts the dependence of the Kolmogorov complexity spectrum, $K_m^C$, on the cell affinity $p$ and logistic parameter $r$ for the coupling parameter $c = 0.02$. This "island" has nearly the same shape and position as the dark shadow square section in Fig. 4 ($r = 4$ and $p = 0.5$), which indicates the asymptotic stability of the system (2). Namely, all around of this section the process of the biochemical substance exchange is highly stochastic because it is not synchronized [28]. In contrast to this, inside of this section the process of the biochemical substance exchange is much less complex indicating absence of stochastic processes and its synchronization.

From panels in Fig. 5 it is seen that for $r = 1$ there exists a large dark gray section that indicates the asymptotic stability, which is expected for the strong coupling. Namely, in that case the dynamical system (2) is not in the chaotic regime and it shows small values of the complexity, since the Kolmogorov complexity $K^C$ is low already for $c = 0.2$ (Fig. 3). When value of $r$ increases, the dark shadow section becomes smaller with larger sections of the asymptotic instability and stability purely depending on nonlinearity of the map $\mathbf{F}$, particularly when approaching to the chaotic regime.

Since the complexity of the processes of exchange of biophysical substances between cells is chaotic ($r = 4$), in Fig. 6 we have plotted stability regions for different values of $c = 0.1, 0.3, 0.5, 0.7, 0.9$ and $p = 0.25, 0.5, 0.75$, where the meaning of axes is the same as in Figs. 4 and 5. We note that the section with the asymptotic stability takes much more place in

($\tilde{x}, \tilde{y}$) space for higher values of $r$. Finally, from all figures it is seen that parameter $r$ influences a much more asymptotic instability than the coupling parameter $c$ does.
Here we consider regions of the asymptotic stability for the equilibrium points in the domain of the map **F** given by (2), through analysis of extreme eigenvalues of the Jacobian, i.e. $\max\{|\lambda_1|, |\lambda_2|\}$ and $\min\{|\lambda_1|, |\lambda_2|\}$ defined by Eqs. (8a)-(8b), in dependence of the equilibrium concentration of substance in the cells ($\tilde{x}, \tilde{y}$) in setting that $r = 4$, $p = 0.5$, $c = 0.02$ (weak coupling) and $c = 0.6$ (strong coupling). The results of computations are depicted in Fig. 7. From Fig. 7a it is seen that for the weak coupling ($c = 0.02$) the minimal values of the $\max\{|\lambda_1|, |\lambda_2|\}$ are reached in the vicinity of the point (0.5, 0.5). Therefore, for such equilibrium point there exists a strongly pronounced asymptotic stability (the green rectangular section). As expected, for the strong coupling ($c = 0.6$) this section of the asymptotic instability is larger (Fig. 7c). Further inspection of Fig. 7, i.e. Fig. (7b) and Fig. (7d) give us information about the asymptotic instability. Namely, for the equilibrium point with cell concentrations close to either higher (towards 1) or lower values (towards 0), values of the $\min\{|\lambda_1|, |\lambda_2|\}$ indicate the existence of a strong instability (green rectangular sections in corners of Fig. 7b). Those sections are smaller in case of stronger coupling (Fig. 7d).

### 3. Inter-cellular biochemical substance exchange in a multi-cell system

*3.1. Stability at an equilibrium state*

In this section we will consider substance exchange in a multi-cell system represented by a ring of coupled cells schematically shown in Fig.8.
In our approach, a cell moves locally in its environment without making long pathways. As a generalization of the two-cell system, according to Mihailovic et al. [6, 38], the dynamics of substance exchange in such a multi-cell system of $N$ cells exchanging biochemical substance, can be represented by the discrete nonlinear time-invariant dynamical system:

$$\mathbf{x}^{(n+1)} = \mathbf{F}(\mathbf{x}^{(n)}) := C\,\Phi(\mathbf{x}^{(n)}) + (I - C)\,Z\,\Psi(\mathbf{x}^{(n)}), \qquad (9)$$

where:
- $x_k^{(n)}$ is the concentration of the substance in $k$-th cell in a discrete time step $n$, $k = 1, 2, ..., N$, $n = 0, 1, 2, ...$, and $\mathbf{x}^{(n)} := \begin{bmatrix} x_1^{(n)} & x_2^{(n)} & ... & x_N^{(n)} \end{bmatrix}^T$ is the appropriate vector,

- $C := diag(c_1, c_2, ..., c_N)$ is the diagonal matrix of the coupling coefficients for each cell,

- $\Phi(x^{(n)}) := diag(\varphi(x_1^{(n)}), \varphi(x_2^{(n)}), ..., \varphi(x_N^{(n)}))$, is the diagonal matrix of intra-cellular behavior modeled by logistic map $\varphi: (0,1) \to (0,1)$, $\varphi(x) := r\,x\,(1-x)$,

- $\Psi(x^{(n)}) := diag\left((x_2^{(n)})^{p_1}, (x_3^{(n)})^{p_2}, ..., (x_N^{(n)})^{p_{N-1}}, (x_1^{(n)})^{p_N}\right)$ is the diagonal matrix of the flow of the substance to each cell, where the all cell's affinities fulfill the constraint $p_1 + p_2 + ... p_N = 1$, and

- $Z \in \{0,1\}^{N \times N}$ is the upper cyclic permutation matrix, i.e., $Z := \begin{bmatrix} e_N & e_1 & e_2 & \ldots & e_{N-1} \end{bmatrix}$, where $e_1, e_2, \ldots, e_N$ are the standard basis vectors of $\mathbb{R}^N$.

It is interesting to note that in this model of multi-cell system for $N$ cells in a ring ($N = 100$ in our simulations) we permit that coupling coefficients can differ from cell to cell. Therefore, two-cell model that we have discussed in the previous section is a special case for $N = 2$ and $c_1 = c_2 = c$.

As in the subsection 2.3, we are interested in the asymptotic stability of the dynamical system (9) at the provided equilibrium point $\tilde{\mathbf{x}} = \mathbf{F}(\tilde{\mathbf{x}})$. Again, since asymptotic stability behavior of nonlinear system (6) at equilibrium is described by the behavior of its linearization around that equilibrium [37], like before, we compute the Jacobian of the map $\mathbf{F}$ at the point $\tilde{\mathbf{x}}$:

$$J_F(\tilde{\mathbf{x}}) := D\mathbf{F}\big|_{\mathbf{x}=\tilde{\mathbf{x}}} = C\,\hat{\Phi}(\tilde{\mathbf{x}}) + (I - C) Z\,\hat{\Psi}(\tilde{\mathbf{x}}), \tag{10}$$

where:
- $\hat{\Phi}(\tilde{\mathbf{x}}) := diag\left(\varphi'(\tilde{x}_1), \varphi'(\tilde{x}_2), \ldots, \varphi'(\tilde{x}_N)\right)$, $\varphi'(x) = r(1 - 2x)$ and
- $\hat{\Psi}(\tilde{\mathbf{x}}) := diag\left(p_1(\tilde{x}_2)^{p_1 - 1}, p_2(\tilde{x}_3)^{p_2 - 1}, \ldots, p_{N-1}(\tilde{x}_N)^{p_{N-1} - 1}, p_N(\tilde{x}_1)^{p_N - 1}\right)$.

In this case, if we use the Theorem from Appendix B, then we obtain that the stability of the discrete dynamical system (9) at the point $\tilde{\mathbf{x}}$, which is determined by the spectral radius $\rho$ of the matrix $J_F(\tilde{\mathbf{x}})$, i.e. $\rho(J_F(\tilde{\mathbf{x}}))$. Thus, the stability of this system is given by

- asymptotically stable if $\rho(J_F(\tilde{\mathbf{x}})) < 1$,

- not asymptotically stable if $\rho\left((J_F(\tilde{\mathbf{x}}))^{-1}\right) > 1$, and

- can be stable or unstable otherwise (due to purely nonlinear effects).

For the specific choice of parameters in the paper, the above criterion gives quite good answer on the question of the stability through computation of the spectra of the Jacobian. Otherwise, finding an analytic expression for the spectral radius $\rho(J_F(\tilde{\mathbf{x}}))$ in dependence on the values of the system parameters is generally impossible. Therefore, in order to have a better insight into the dependence of stability on the parameter changes, we will return to matrix infinity norm. Namely, it is well known that

$$\rho(J_F(\tilde{\mathbf{x}})) \leq \|J_F(\tilde{\mathbf{x}})\|_\infty = \max_{i=1,2,\ldots,N} \sum_{j=1}^{N} \left|[J_F(\tilde{\mathbf{x}})]_{i,j}\right|. \tag{11}$$

Having this in mind, we obtain only the necessary conditions for the asymptotic stability of the system (6) by computing the infinity norm and analyzing the following inequality:

$$\|J_F(\tilde{\mathbf{x}})\|_\infty = \max_{i=1,2,\ldots,N} \left((1 - c_i)\,r\,(1 - 2\tilde{x}_i) + c_i\,p_i(\tilde{x}_{i+1})^{p_i - 1}\right) < 1. \tag{12}$$

At this point, let us note that the specific choice of infinity norm, in order to bound the spectral radius of a matrix, is not generally the best possible choice. In fact, this is equivalent to approximations of matrix eigenvalues by Geršgorin's circles (see [39]). But, due to its analytic simplicity it is a good choice in many cases, as a starting point. The better estimations of stability region could be done by the use of other matrix norms or localization areas [40-43], which will be in the focus of our future work regarding this subject.

*3.2. Domains for equilibrium points that permit stability for every coupling*

We start by observing that inequality (12) holds if and only if for every $i \in \{1, 2, ..., N\}$

$$(1-c_i) r |1-2\tilde{x}_i| + c_i p_i (\tilde{x}_{i+1})^{p_i - 1} < 1. \tag{13}$$

But, since the left hand side of (13) is a convex combination of $r|1-2\tilde{x}_i|$ and $p_i(\tilde{x}_{i+1})^{p_i-1}$ with parameter $0 < c_i < 1$, we conclude that $r|1-2\tilde{x}_i| < 1$ and $p_i(\tilde{x}_{i+1})^{p_i-1} < 1$ together imply (12) independently of the values of $0 < c_i < 1$, $i \in \{1, 2, ..., N\}$.

Thus, since $r|1-2\tilde{x}_i| < 1$ and $p_i(\tilde{x}_{i+1})^{p_i-1} < 1$ is equivalent to:

$$\tilde{x}_i \in \left( \frac{r-1}{2r}, \frac{r+1}{2r} \right) \quad \text{and} \quad (p_{i-1})^{\frac{1}{1-p_{i-1}}} < \tilde{x}_i, \tag{14}$$

we obtain the region in the $N$-dimensional space of substance concentrations in the coupled ring of cells

$$\mathbf{S} := \left\{ \tilde{\mathbf{x}} \in (0,1)^N : \max \left\{ \frac{r-1}{2r}, (p_{i-1})^{\frac{1}{1-p_{i-1}}} \right\} < \tilde{x}_i < \frac{r+1}{2r} \right\}, \tag{15}$$

such that for every coupling ($0 < c_i < 1, i \in \{1, 2, ..., N\}$), if an equilibrium point $\tilde{\mathbf{x}} \in \mathbf{S}$, then this equilibrium point is asymptotically stable. In other words, independently of the coupling parameters of the individual cells in the multi-cell system, region $\mathbf{S}$, which always contains a hypercube around the central point the domain $[0.5 \ 0.5 \ ... \ 0.5]^T$ with the edge of the length $0.5r^{-1}$, is the place where asymptotic stability of the point of equilibrium is always assured. Looking at Fig. 6, we can see the effect of the same behavior in the two-cell system, i.e., while coupling parameter was changing the square region $(0.375, 0.625) \times (0.375, 0.625)$ somehow remains included in the stability region. However, using an infinity norm bound, in fact, we can explain existence of this region even for system of $N$ cells coupled in a ring formation. Namely, if we assume that the values of the logistic parameter $r$ are ranged from $3.785$ to $4$, then, after a closer look at (15), it is seen that, independently on the coupling parameters and the cell affinities, the following inclusions hold

$$(0.375, 0.625)^N \subseteq \mathbf{S} \subseteq (0.3679, 0.6321)^N. \tag{16}$$

Namely, if an equilibrium concentration in each cell is between 0.375 and 0.625 then asymptotic stability of (9) is assured without any additional constraints on the parameters ($c_i, r, p_i$), i.e., for each $c_1, c_2, ..., c_N \in (0,1)$, $p_1, p_2, ..., p_N \in (0,1)$ and $r \in (3.785, 4)$. In addition, we consider a case when $p_i \leq 0.8, i \in \{1, 2, ..., N\}$. Similarly, as in above analysis of (15) we conclude that this constraint extends the stability interval for individual cell's equilibrium concentration to

$$(0.375, 0.625)^N \subseteq \mathbf{S} \subseteq (0.3333, 0.6667)^N. \qquad (17)$$

Here, the constraint $p_i \leq 0.8, i \in \{1, 2, ..., N\}$, encounters practically in majority of the cases of the multi-cell systems since $p_i + p_2 + ... + p_N = 1$.

Finally, we underline a conclusion, which can be derived from the above discussion. Namely, the asymptotic stability of the dynamic system represented by a ring of coupled cells given by (9) is assured without any additional constraints for $r \in (3.785, 4)$, i.e. the interval, which includes the size of the $r$ interval for the "island" of the low complexity for two cell system exchanging the biochemical substance (Fig. 2a). This is interesting since in this interval the coupled maps (4) show a chaotic behavior. It means that in those conditions, there exists space of parameters ($c_i, r, p_i$), for which the process of biochemical substance exchange in a coupled ring of cells is stable.

## 5. Concluding remarks

In this paper we have considered the complexity and asymptotic stability in the process of biochemical substance exchange in a coupled ring of cells. To model this process we have used a coupled map, which includes a formalization of the cell affinity in the interaction with other cells. We have given a short outline of the inter-cellular biochemical substance exchange model. We have considered: (i) the complexity (by introducing the Kolmogorov complexity spectrum and the Kolmogorov complexity spectrum highest value and) and (ii) asymptotic stability (using an eigenvalue-based method) of the dynamical system for two cells. Then we have explored how stability of the equilibrium of the biochemical substance exchange in a multi-cell system ($N = 100$) is influenced by the coupling parameter, cell affinity and environmental factor as master parameters of the model. We have found that even in highly chaotic conditions there exists space of these parameters, for which the process of biochemical substance exchange in a coupled ring of cells is stable.

**Acknowledgements**

This paper was realized as a part of the projects "Studying climate change and its influence on the environment: impacts, adaptation and mitigation" (43007) and "Numerical linear algebra and discrete structures" (174019) financed by the Ministry of Education and Science of the Republic of Serbia within the framework of integrated and interdisciplinary research for the period 2011-2014, and as a part of the project "Applied linear and combinatorial algebra"

(2675) financed by the Provincial Secretariat for Science and Technological Development of AP Vojvodina within the project cycle 2011-2014.

# Appendix A

*Description of the Lempel- Ziv algorithm for computing the Kolmogorov complexity*

Lempel- Ziv algorithm [35] for computing the Kolmogorov complexity of a time series $\{x_i\}$, $i = 1,2,3,4,...,N$, can be summarized as follows. *Step A:* Encode the time series by constructing a sequence $S$ consisting of the characters 0 and 1 written as $\{s(i)\}$, $i=1,2,3,4,…,N$, according to the rule

$$s(i) = \begin{cases} 0 & x_i < x_* \\ 1 & x_i \geq x_* \end{cases}. \tag{A1}$$

Here $x_*$ is a threshold that should be properly chosen. The mean value of the time series has often been used as the threshold [44]. Depending on the application, other encoding schemes are also available [45]. *Step B:* Calculate the complexity counter $C(N)$, which is defined as the minimum number of distinct patterns contained in a given character sequence [46]; $c(N)$ is a function of the length of the sequence $N$. The value of $c(N)$ is approaching an ultimate value $b(N)$ as $N$ approaching infinite, i.e.

$$c(N) = O(b(N)), \; b(N) = \frac{N}{\log_2 N}. \tag{A2}$$

*Step C*: Calculate the normalized complexity measure $C_k(N)$, which is defined as

$$C_k(N) = \frac{c(N)}{b(N)} = c(N) \frac{\log_2 N}{N}. \tag{A3}$$

The $C_k(N)$ is a parameter to represent the information quantity contained in a time series, and it should be a 0 for a periodic or regular time series or 1 for a random time series, if $N$ is large enough. For a non-linear time series, $C_k(N)$ is to be between 0 and 1.

# Appendix B

*Stability of the discrete nonlinear dynamical system through linearization*

Linearization technique is generally used in obtaining asymptotic stability of discrete dynamical systems, [37]. Here we give one such Theorem and its proof that clarify the analysis in the previous sections.

**Theorem:** *Given a nonlinear continuously Fréchet differentiable map* $\mathbf{F} : (0,1)^N \to (0,1)^N$ *and the fixed point* $\tilde{\mathbf{x}} = \mathbf{F}(\tilde{\mathbf{x}})$ *of the discrete nonlinear dynamical system*

$$\mathbf{x}^{(n+1)} = \mathbf{F}(\mathbf{x}^{(n)}), \quad n = 0, 1, 2, \ldots \tag{B1}$$

let $J_F(\tilde{\mathbf{x}}) := D\mathbf{F}\big|_{\mathbf{x}=\tilde{\mathbf{x}}}$ denote the Jacobian matrix at the point $\tilde{\mathbf{x}}$. Then, the following two implications hold:

(a) if $\rho(J_F(\tilde{\mathbf{x}})) < 1$, then the system (B1) is asymptotically stable at the point $\tilde{\mathbf{x}}$,

(b) if all eigenvalues of $J_F(\tilde{\mathbf{x}})$ are larger in absolute value than one, i.e., $J_F(\tilde{\mathbf{x}})$ is invertible and $\rho\big((J_F(\tilde{\mathbf{x}}))^{-1}\big) > 1$, then the system (B1) is not asymptotically stable at the point $\tilde{\mathbf{x}}$.

**Proof:** First, supposing that the map $\mathbf{F}$ is continuously Fréchet differentiable on the domain, we obtain that the Jacobian matrix exists at each point of the domain, that it is a continuous map, and that the Fréchet derivative is given by:

$$\mathbf{F}'_{\mathbf{x}}(h) = J_F(\mathbf{x})h, \quad h \in (0,1)^N.$$

Now, assume that $\rho(J_F(\tilde{\mathbf{x}})) < 1$, then, see [39], there exists a matrix norm $\|\ \|$ on the space $(0,1)^{N,N} \subset \mathbb{R}^{N,N}$ induced by some vector norm $\|\ \|_v$, such that $\rho(J_F(\tilde{\mathbf{x}})) < \|J_F(\tilde{\mathbf{x}})\| < 1$. Therefore, we obtain $\|\mathbf{F}'_{\tilde{\mathbf{x}}}\| = \|J_F(\tilde{\mathbf{x}})\| = \sup\{\|J_F(\tilde{\mathbf{x}})h\|_v : \|h\|_v = 1\} < 1$.

Now, since every norm is a continuous, we have that the map $\mathbf{x} \to \|\mathbf{F}'_{\mathbf{x}}\|$ is a continuous map, too. Therefore, there exist $\varepsilon > 0$ and $M \in [0,1)$ such that:

$$\|\mathbf{F}'_{\mathbf{x}}\| < M, \text{ for all } \mathbf{x} \in \mathbf{B}(\tilde{\mathbf{x}}, \varepsilon) := \{\mathbf{y} \in (0,1)^N : \|\tilde{\mathbf{x}} - \mathbf{y}\| < \varepsilon\}.$$

Further, using the mean value theorem we have that $\|\mathbf{F}(\mathbf{x}) - \mathbf{F}(\mathbf{y})\|_v \leq \|\mathbf{F}'_{\mathbf{z}}\| \|\mathbf{x} - \mathbf{y}\|_v$, for all $\mathbf{x}, \mathbf{y} \in \mathbf{B}(\tilde{\mathbf{x}}, \varepsilon)$, and some $\mathbf{z}$ which is a convex combination of $\mathbf{x}$ and $\mathbf{y}$, and, therefore $\mathbf{z} \in \mathbf{B}(\tilde{\mathbf{x}}, \varepsilon)$. But, this implies that

$$\|\mathbf{F}(\mathbf{x}) - \tilde{\mathbf{x}}\|_v = \|\mathbf{F}(\mathbf{x}) - \mathbf{F}(\tilde{\mathbf{x}})\|_v \leq M \|\mathbf{x} - \tilde{\mathbf{x}}\|_v < \varepsilon \text{ for all } \mathbf{x} \in \mathbf{B}(\tilde{\mathbf{x}}, \varepsilon), \text{ i.e.,}$$

$$\mathbf{F}\big[\mathbf{B}(\tilde{\mathbf{x}}, \varepsilon)\big] \subseteq \mathbf{B}(\tilde{\mathbf{x}}, \varepsilon).$$

Therefore, if the starting point $\mathbf{x}^{(0)}$ of the discrete dynamical system (B1) is in the $\varepsilon$-neighborhood of the fixed point $\tilde{\mathbf{x}}$, then for $\mathbf{x}^{(n)}$, the state of (B1) at $n$-th discrete time step, we have that

$$\|\mathbf{x}^{(n)} - \tilde{\mathbf{x}}\|_v = \|\mathbf{F}(\mathbf{x}^{(n-1)}) - \tilde{\mathbf{x}}\|_v \leq M\|\mathbf{x}^{(n-1)} - \tilde{\mathbf{x}}\|_v = \ldots \leq M^n\|\mathbf{x}^{(0)} - \tilde{\mathbf{x}}\|_v < M^n \varepsilon.$$

and, as a consequence, when $n \to \infty$, $\mathbf{x}^{(n)} \to \tilde{\mathbf{x}}$, and the asymptotic stability is obtained.

In order to prove item (b), we observe that if $\rho\left(\left[J_F(\tilde{\mathbf{x}})\right]^{-1}\right) < 1$, according to the same results of [39], then there exist vector norm $\|\ \|_v$ and its induced matrix norm $\|\ \|$, such that $\rho\left(\left[J_F(\tilde{\mathbf{x}})\right]^{-1}\right) < \left\|\left[J_F(\tilde{\mathbf{x}})\right]^{-1}\right\| = \left\|\left[\mathbf{F'}_{\tilde{\mathbf{x}}}\right]^{-1}\right\| < 1$. Thus, we have obtained that $\left\|\left[\mathbf{F'}_{\tilde{\mathbf{x}}}\right]^{-1}\right\|^{-1} > 1$.

But, then, for every $\mathbf{x}, \mathbf{y} \in \mathbf{B}(\tilde{\mathbf{x}}, \varepsilon)$,

$$\|\mathbf{x} - \mathbf{y}\|_v = \left\|\left[\mathbf{F'}_{\tilde{\mathbf{x}}}\right]^{-1}\left(\mathbf{F'}_{\tilde{\mathbf{x}}}(\mathbf{x})\right) - \left[\mathbf{F'}_{\tilde{\mathbf{x}}}\right]^{-1}\left(\mathbf{F'}_{\tilde{\mathbf{x}}}(\mathbf{y})\right)\right\|_v \le \left\|\left[\mathbf{F'}_{\tilde{\mathbf{x}}}\right]^{-1}\right\| \left\|\mathbf{F'}_{\tilde{\mathbf{x}}}(\mathbf{x} - \mathbf{y})\right\|_v,$$

which implies that $\left\|\mathbf{F'}_{\tilde{\mathbf{x}}}(\mathbf{x} - \mathbf{y})\right\|_v \ge q \|\mathbf{x} - \mathbf{y}\|_v$, for $q = \left\|\left[\mathbf{F'}_{\tilde{\mathbf{x}}}\right]^{-1}\right\|^{-1} > 1$.

Since the map $\mathbf{F}$ is Fréchet differentiable, we have that

$$\mathbf{F}(\mathbf{x}) - \mathbf{F}(\tilde{\mathbf{x}}) = \mathbf{F'}_{\tilde{\mathbf{x}}}(\mathbf{x} - \tilde{\mathbf{x}}) + o\left(\|\mathbf{x} - \tilde{\mathbf{x}}\|_v\right) \text{ for all } \mathbf{x} \in (0,1)^N,$$

and, thus, we can choose $\delta > 0$ such that $\hat{q} = q - \delta > 1$, and $\varepsilon > 0$ such that for every $\mathbf{x} \in \mathbf{B}(\tilde{\mathbf{x}}, \varepsilon) \setminus \{\tilde{\mathbf{x}}\}$, $\left\|o\left(\|\mathbf{x} - \tilde{\mathbf{x}}\|_v\right)\right\|_v \le \delta \|\mathbf{x} - \tilde{\mathbf{x}}\|_v$, and obtain that

$$\|\mathbf{F}(\mathbf{x}) - \tilde{\mathbf{x}}\|_v = \|\mathbf{F}(\mathbf{x}) - \mathbf{F}(\tilde{\mathbf{x}})\|_v \ge \left\|\mathbf{F'}_{\tilde{\mathbf{x}}}(\mathbf{x} - \tilde{\mathbf{x}})\right\|_v - \left\|o\left(\|\mathbf{x} - \tilde{\mathbf{x}}\|_v\right)\right\|_v \ge \hat{q}\|\mathbf{x} - \tilde{\mathbf{x}}\|_v - \delta\|\mathbf{x} - \tilde{\mathbf{x}}\|_v = q\|\mathbf{x} - \tilde{\mathbf{x}}\|_v,$$

for all $\mathbf{x} \in \mathbf{B}(\tilde{\mathbf{x}}, \varepsilon) \setminus \{\tilde{\mathbf{x}}\}$.

Now, let $\mathbf{x}^{(0)}$ be an arbitrary starting point of the discrete dynamical system (B1) is in the $\varepsilon$-neighborhood of the fixed point $\tilde{\mathbf{x}}$, such that all the states of (B1) up to the $n-1$–th discrete time step all belong to that $\varepsilon$-neighborhood of $\tilde{\mathbf{x}}$. Then, we have that

$$\left\|\mathbf{x}^{(n)} - \tilde{\mathbf{x}}\right\|_v = \left\|\mathbf{F}\left(\mathbf{x}^{(n-1)}\right) - \tilde{\mathbf{x}}\right\|_v \ge q\left\|\mathbf{x}^{(n-1)} - \tilde{\mathbf{x}}\right\|_v \ge \ldots \ge q^n \left\|\mathbf{x}^{(0)} - \tilde{\mathbf{x}}\right\|_v,$$

which implies that there must exist (sufficiently large) discrete time step $n$ such that $\mathbf{x}^{(n)} \notin \mathbf{B}(\tilde{\mathbf{x}}, \varepsilon)$. Therefore, in the chosen $\varepsilon$-neighborhood of the fixed point $\tilde{\mathbf{x}}$ there cannot exist a sequence of states of (B1) that will converge to $\tilde{\mathbf{x}}$, i.e., (B1) cannot be asymptotically stable at the point $\tilde{\mathbf{x}}$.

**Figures and Captions**

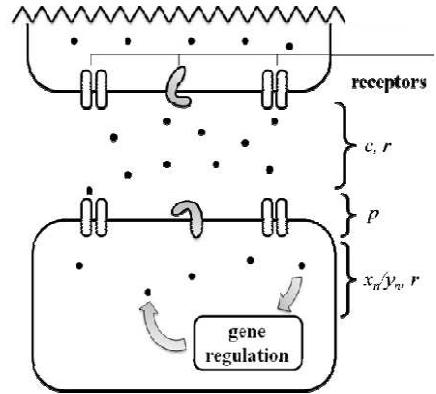

**Fig. 1 Schematic representation of intercellular communication [reprinted from 28]. Here, $c$ represents concentration of signaling molecules in intercellular environment coupled with intensity of response they can provoke while $r$ includes collective influence of environment factors which can interfere with the process of communication. $x_n$ and $y_n$ represent concentration of signaling molecules in cells environment, while $p$ represents the cell affinity.**

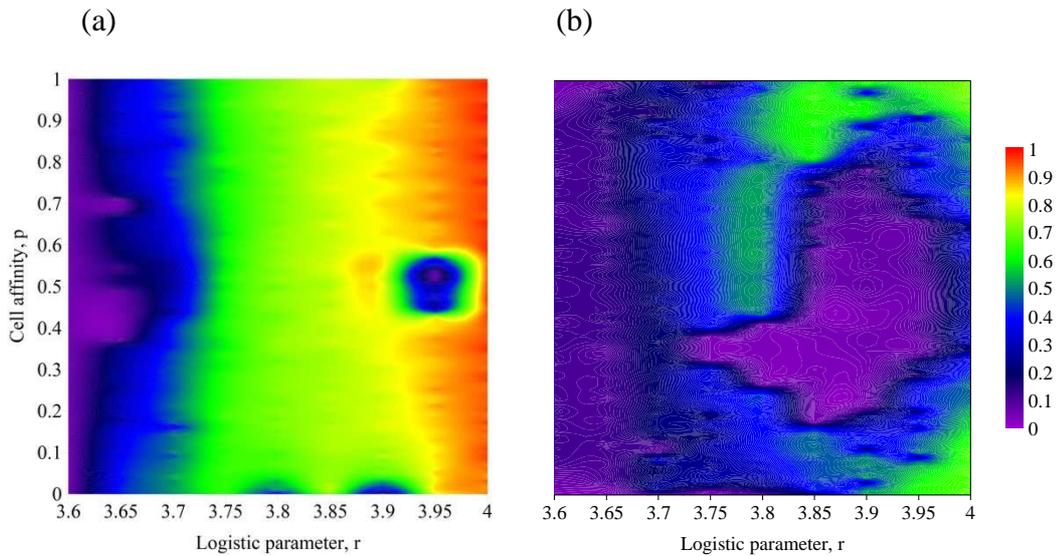

**Fig. 2 The dependence of the Kolmogorov complexity spectrum highest value, $K_m^C$ on the cell affinity $p$ and logistic parameter $r$ for the coupling parameter $c$: (a) 0.02 and (b) 0.2.**

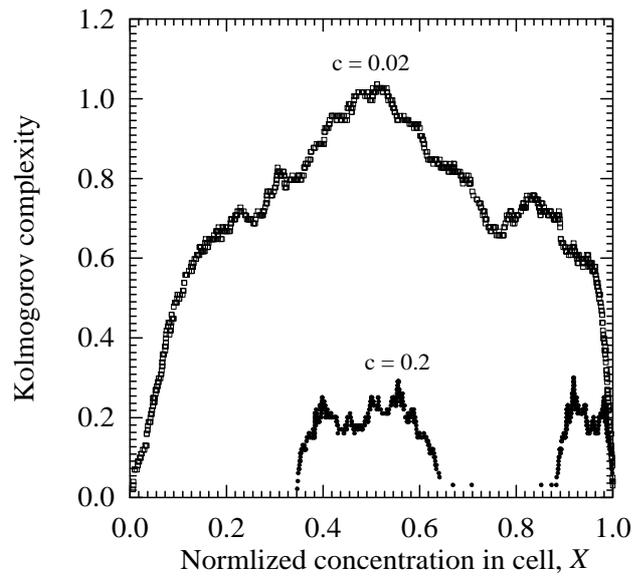

Fig. 3 Kolmogorov complexity spectrum of the process of the biochemical substance exchange between two cells (Eq. (4a)-(4b)) as a function of normalized concentration in cell, $X$ for the coupling parameter $c$ : 0.02 and 0.2.

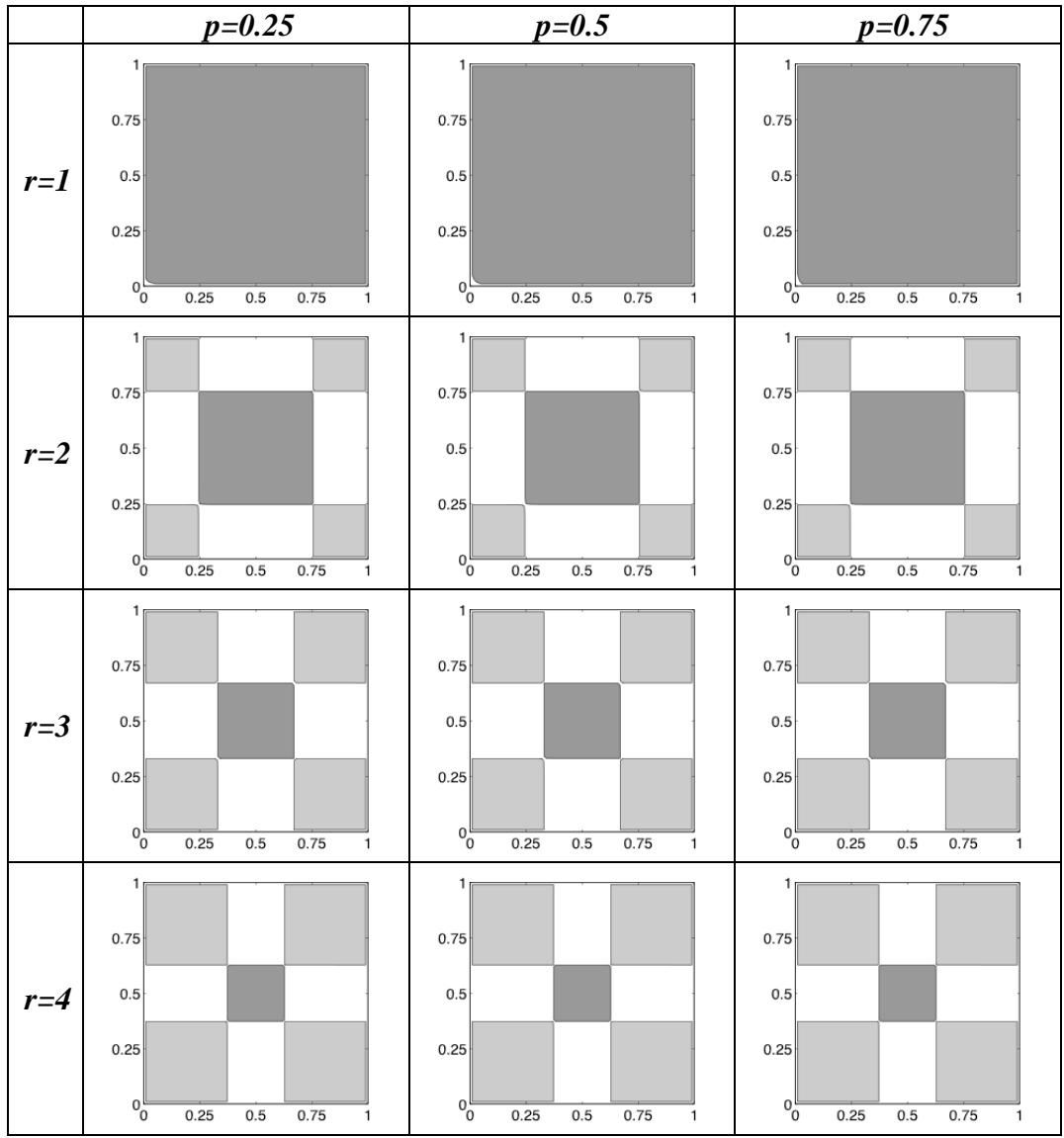

**Fig. 4** Regions of the asymptotic stability for the equilibrium points (Eqs. (8)), in the domain of the map F given by (2), for weak coupling $c = 0.02$. The values of the equilibrium concentration of substance in the first cell ($\tilde{x}$) and second one ($\tilde{y}$) are given on the horizontal and vertical axes, respectively. Domains of asymptotic stability and instability are indicated by dark and light gray areas, respectively. White area indicates that stability purely depends on nonlinearity of the map F around the point of equilibrium $\tilde{v}$.

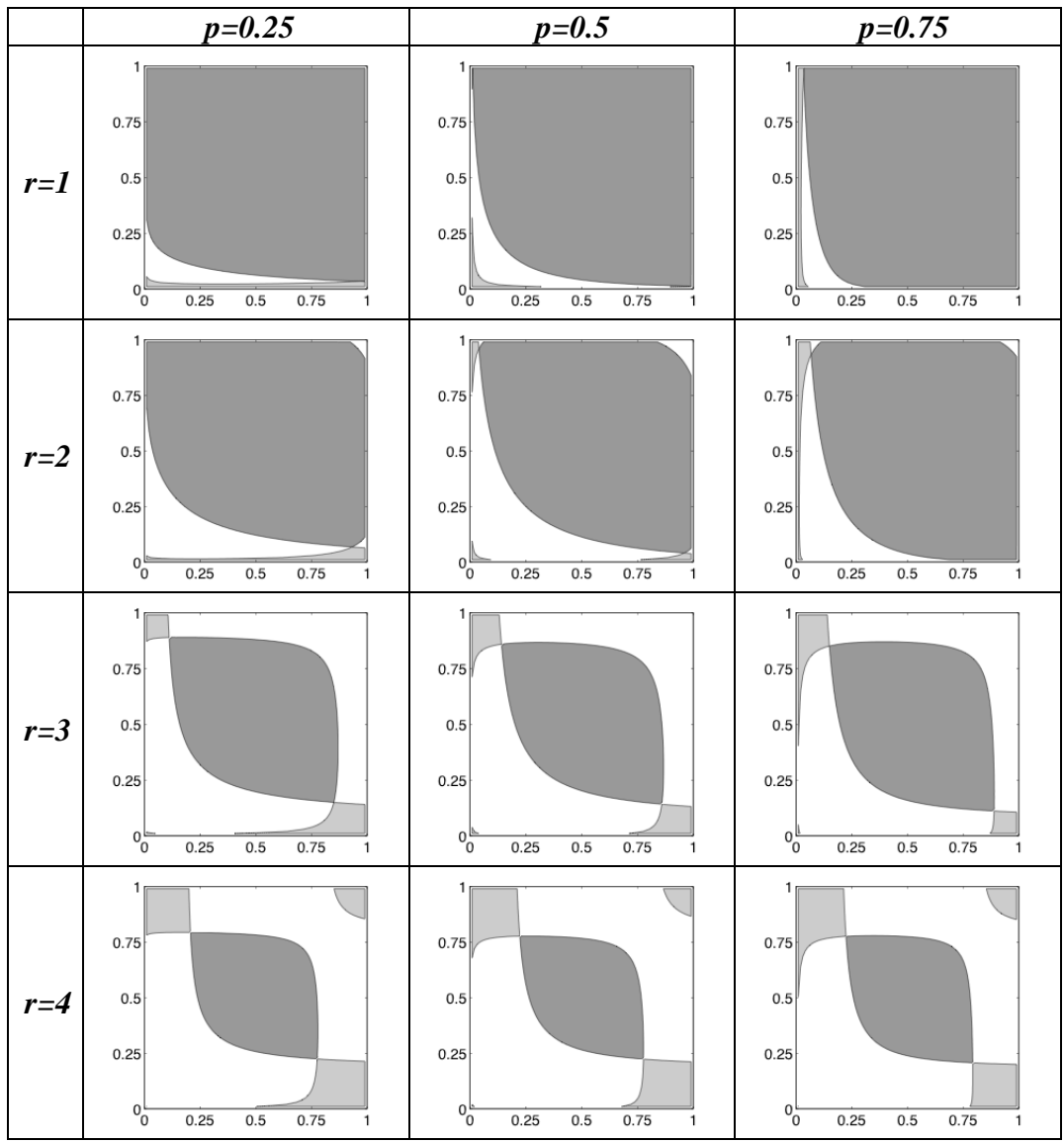

**Fig. 5** The same as Fig. 4 but for strong coupling *c=0.6*.

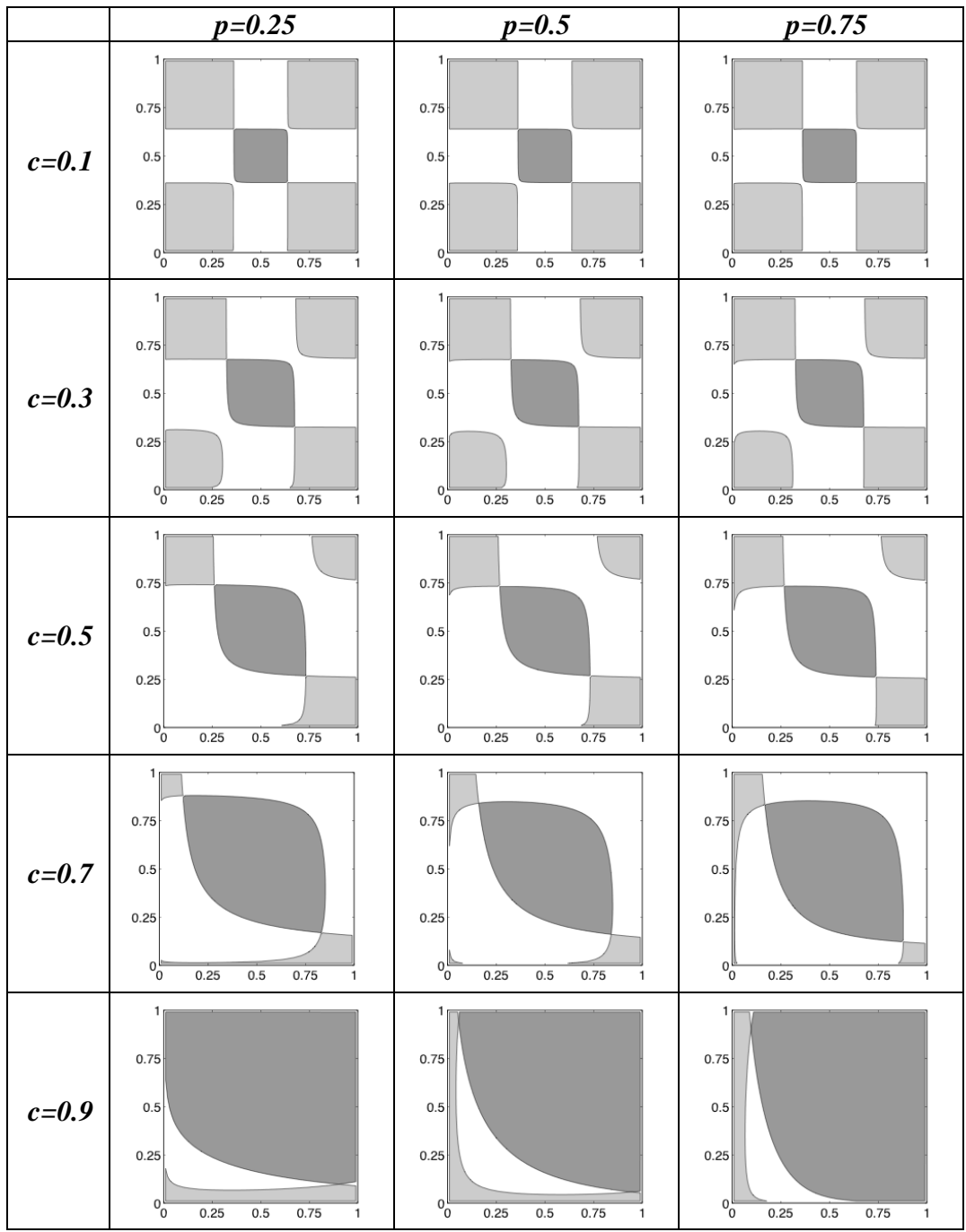

**Fig. 6** Regions of the asymptotic stability for the equilibrium points (Eqs. (8)), in the domain of the map F given by (2), for $r = 4$ and different values of coupling parameter $c$.

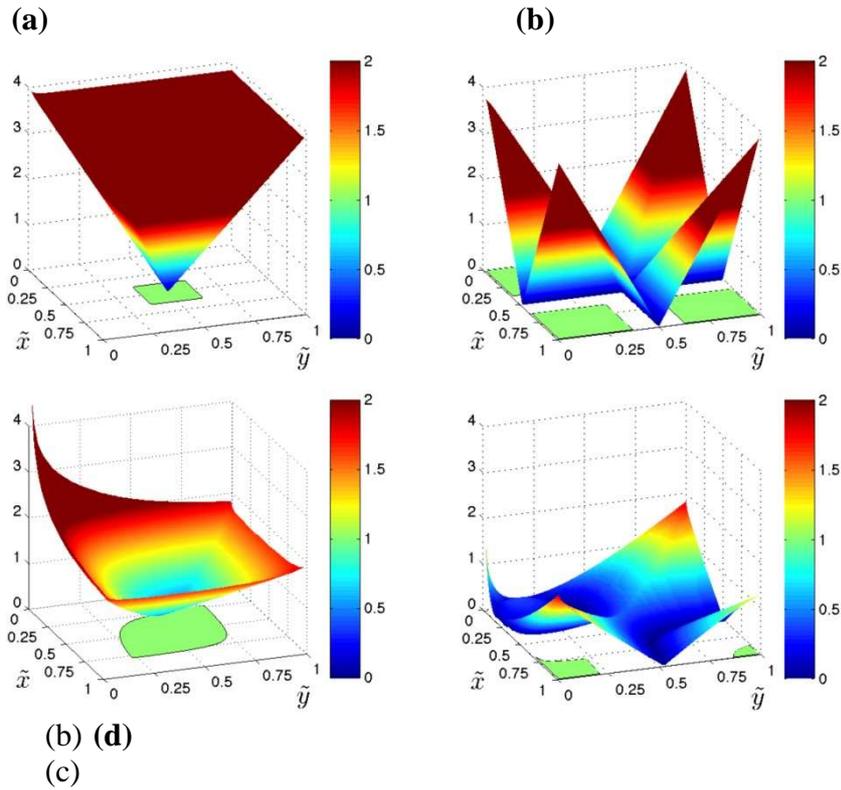

**Fig. 7** 3D diagram of the extreme eigenvalues of the Jacobian, $\max\{|\lambda_1|,|\lambda_2|\}$ ((a) and (c)) and $\min\{|\lambda_1|,|\lambda_2|\}$ ((b) and (d)), against equilibrium concentration of substance in the cells ($\tilde{x}$, $\tilde{y}$) in setting: $r = 4$, $p = 0.5$, $c = 0.02$ (weak coupling in (a) and (b)) and $c = 0.6$ (strong coupling in (c) and (d)). Sections in ($\tilde{x}$, $\tilde{y}$) plane are due to the regions with stability and instability.

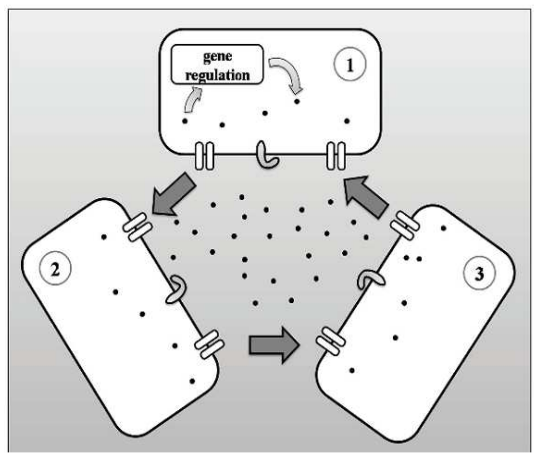

**Fig 8.** Schematic diagram of a model of substance exchange in a system represented by a ring of coupled cells [reprinted from 6]